\pdfoutput=1
\RequirePackage{ifpdf}
\ifpdf % We~are running pdfTeX in pdf mode
\documentclass[pdftex]{sigma}
\else
\documentclass{sigma}
\fi

\numberwithin{equation}{section}

\newtheorem{Theorem}{Theorem}[section]
\newtheorem*{Theorem*}{Theorem}
\newtheorem{Corollary}[Theorem]{Corollary}
\newtheorem{Conjecture}[Theorem]{Conjecture}
\newtheorem{Lemma}[Theorem]{Lemma}
\newtheorem{Proposition}[Theorem]{Proposition}

\theoremstyle{definition}
\newtheorem{Definition}[Theorem]{Definition}

\newtheorem{Example}[Theorem]{Example}
\newtheorem{Remark}[Theorem]{Remark}

\usepackage{bm}

\let\frak\mathfrak

\global\long\def\<{\langle}

\global\long\def\>{\rangle}

\global\long\def\l{\lambda}

\global\long\def\g{\frak{g}}

\global\long\def\hg{\hat{\frak{g}} }
\global\long\def\hb{\hat{\frak{b}} }

\global\long\def\h{\frak{h}}

\global\long\def\n{\frak{n}}
\global\long\def\Z{\mathbb{Z}}

\global\long\def\Zp{{\mathbb Z}_{\ge0} }
\global\long\def\N{\mathbb{N}}
\global\long\def\C{\mathbb{C}}

\global\long\def\sl{\mathfrak{sl}}

\begin{document}

\allowdisplaybreaks

\newcommand{\arXivNumber}{2512.02718}

\renewcommand{\thefootnote}{}

\renewcommand{\PaperNumber}{052}

\FirstPageHeading

\ShortArticleName{Irreducibility of Certain $\widehat{\mathfrak sl}_2$-Modules of Wakimoto Type}

\ArticleName{Irreducibility of Certain $\boldsymbol{\widehat{\mathfrak{sl}}_2}$-Modules\\ of Wakimoto Type\footnote{This paper is a~contribution to the Special Issue on Recent Advances in Vertex Operator Algebras in honor of James Lepowsky. The~full collection is available at \href{https://sigma-journal.com/Lepowsky.html}{https://sigma-journal.com/Lepowsky.html}}}

\Author{Dra\v{z}en ADAMOVI\'C and Veronika PEDI\'C TOMI\'C}

\AuthorNameForHeading{D.~Adamovi\'c and V.~Pedi\'c Tomi\'c}

\Address{Department of Mathematics, Faculty of Science, University of Zagreb,\\
Bijeni\v{c}ka 30, Zagreb, Croatia}
\Email{\mail{adamovic@math.hr}, \mail{vpedic@math.hr}}
\URLaddress{\url{https://web.math.pmf.unizg.hr/~adamovic/},\newline
\hspace*{10.5mm}\url{https://web.math.pmf.unizg.hr/~vpedic/}}

\ArticleDates{Received December 03, 2025, in final form May 11, 2026; Published online May 21, 2026}

\Abstract{We investigate the irreducible smooth $\widehat{\mathfrak{sl}}_{2}$-modules recently constructed in [\textit{Adv. Math.} \textbf{481} (2025), 110559, 34~pages, arXiv:2404.03855], and demonstrate that these modules admit a Wakimoto-type realization at both critical and non-critical levels. In the critical level case, we identify simple quotients of these modules with the Wakimoto modules whose irreducibility was already established by Adamovi\'c. We~also generalize some Wakimoto modules constructed in [\textit{Adv. Math.} \textbf{289} (2016), 438--479, arXiv:1409.5354] and identify them as generalized Whittaker modules.}


\Classification{17B69; 17B20; 17B67; 17B68}

\begin{flushright}
\begin{minipage}{60mm}
\it Dedicated to James Lepowsky\\
 on the occasion of his 80th birthday
\end{minipage}
\end{flushright}

\renewcommand{\thefootnote}{\arabic{footnote}}
\setcounter{footnote}{0}

\section{Introduction}

The Wakimoto modules are representations of affine Kac--Moody Lie algebras.
They were first introduced by M.~Wakimoto \cite{Wak} in the case of $\widehat{\mathfrak{sl}}_{2}$
and later generalized by E.~Frenkel and B.~Feigin \cite{FF,FF2,FF1}.
A vertex-algebraic approach to Wakimoto modules was described in \cite{efren}.
 Wakimoto modules were also used for a construction of positive energy modules (cf. \cite{FK}) and for the determination of the associated varieties of affine vertex algebras (cf. \cite{AFK}). A version of Wakimoto modules for studying imaginary highest weight modules (which are not smooth) was introduced in \cite{CF, CF2}.

Let $V^{\kappa}(\g)$ denote the universal affine vertex algebra associated with a simple Lie algebra $\g$ at level $\kappa$.
Assume that $\g = {\n}_{-} \oplus {\h} \oplus {\n}_{+}$ is the usual triangular decomposition of $\g$.
Then there exists an injective homomorphism of vertex algebras, known as the Frenkel--Feigin homomorphism,
\[
 \Phi^{\kappa} \colon \ V^{\kappa}(\g) \longrightarrow W_{\g} \otimes \pi_{\h}^{\kappa + h^{\vee}},
\]
where $W_{\g}$ is the Weyl vertex algebra associated with $\g$, and
$\pi_{\h}^{\kappa + h^{\vee}}$ is the Heisenberg vertex algebra of level $\kappa + h^{\vee}$ associated with $\h$
(see \cite{efren} for details).

In particular, $W_{\g}$ is the tensor product of $\dim \n_{-}$ copies of the rank-one Weyl vertex al\-gebra~$W$,
while \smash{$\pi_{\h}^{\kappa + h^{\vee}}$} a Heisenberg vertex algebra of rank equal to $\dim \h$.

Suppose that $M_1$ is a $W_{\g}$-module and $N_1$ a
\smash{$\pi_{\h}^{\kappa + h^{\vee}}$}-module.
Then, using the Frenkel--Feigin homomorphism, the tensor product $M_1 \otimes N_1$ naturally acquires the structure of a weak $V^{\kappa}(\g)$-module.
This module is called the \emph{Wakimoto module}, since in the case $\mathfrak{sl}_2$ this construction was originally discovered by M.~Wakimoto (cf.~\cite{Wak}).
If $N_1$ is taken to be a highest weight module for the Heisenberg vertex algebra, then Wakimoto modules provide a framework for studying highest weight modules for affine Lie algebras.
For generic highest weights, Wakimoto modules are isomorphic to irreducible highest weight modules of affine Lie algebras.
The irreducibility of certain Wakimoto modules gave a natural proof of the Kac--Kazhdan conjecture \cite{KK} concerning the characters of irreducible representations at the critical level
(see \cite[Section~4]{efren} for details).

The connection with Kac--Kazhdan conjecture, motivated the investigation of irreducibility Wakimoto modules of non-highest weight type in \cite{A-07, A-14, AJM}. These results are obtained for the affine vertex algebra $V^{\kappa}(\g)$ at the critical level. The natural step is to consider analogs of these modules of non-critical level. In this paper, we propose one generalization.

We may consider the special case $M_1 = W_{\mathfrak{sl}_2}$ and $N_1$ an arbitrary $\pi_{\h}^{\kappa + h^{\vee}}$-module.
It is then natural to ask under what conditions the tensor product $M_1 \otimes N_1$ becomes an irreducible $V^{\kappa}(\g)$-module.

In this paper, we first study the case of non-critical levels, and we propose the following conjecture.

\begin{Conjecture} \label{slutnja}
Assume that $\kappa$ is non-critical and that
$N_1(\bm{\lambda})$ is a Whittaker, non-highest-weight module for the Heisenberg vertex algebra \smash{$\pi^{\kappa + h^{\vee}}$}.
Then $W_{\g} \otimes N_1(\bm{\lambda})$ is an irreducible $V^{\kappa}(\g)$-module.
\end{Conjecture}

In Theorem~\ref{prvi-teorem}, we prove Conjecture~\ref{slutnja} in the case $\g = \mathfrak{sl}_2$:
\begin{itemize}\itemsep=0pt
 \item If $N_1(\bm{\lambda})$ is a Whittaker module for the Heisenberg vertex algebra \smash{$\pi^{\kappa + h^{\vee}}$},
 then $W \otimes N_1(\bm{\lambda})$ is an irreducible $V^{\kappa}(\mathfrak{sl}_2)$-module.
\end{itemize}

We further show that the above Wakimoto modules for $\g = \mathfrak{sl}_2$ are isomorphic to the irreducible modules
\smash{$\widehat{M}(\varphi)$} recently studied by V.~Futorny, X.~Guo, Y.~Xue, and K.~Zhao in \cite{FGXZ} as certain universal modules
(cf.~Section~\ref{FGXZ-approach} for a precise definition).

A special interest arises in the case of the critical level, when $\pi_{\h}^{\kappa + h^{\vee}}$ becomes a commutative vertex algebra.
In \cite{A-07,A-14}, the first-named author considered Wakimoto modules for $V^{\kappa}(\mathfrak{sl}_2)$ at the critical level
and completely solved the question of their irreducibility.
We recall this result in Section~\ref{ired-critical-sl(2)}.
At the critical level, we identify the quotients $\overline{M}(\varphi, \theta)$ of the modules \smash{$\widehat{M}(\varphi)$}
as Wakimoto modules parameterized by $\chi \in \mathbb{C}((z))$ (cf.~Theorem~\ref{krit-kvocijenti})
where
the center of $V^{-2} (\sl_2)$ acts as $ T(z) = \theta(z) = \sum _{n \in {\Z}} \theta_n z^{-n-2} \in \C((z)).$

Our results slightly generalize those of \cite{FGXZ}, showing that in certain cases the modules $\overline{M}(\varphi, \theta)$ are reducible
(cf.~Proposition~\ref{reducibilni}).
The next natural step is to consider modules of the form $M_1 \otimes N_1$, where $M_1$ is a Whittaker module for $W$ (cf.~\cite{ALPY}), and $N_1$ is a Whittaker \smash{$\pi^{\kappa +h^{\vee}}$}-module. In Section~\ref{gen-whitt-sect}, we identify such modules as generalized Whittaker modules. At the critical level, we establish new results concerning their irreducibility. A more detailed analysis will be presented in our forthcoming paper~\cite{AP-2026}.

\subsection*{Notation}
\begin{itemize}\itemsep=0pt
 \item $\hg = \g\otimes \bigl[t, t^{-1}\bigr]\oplus \C K$.
 \item $\h\subset \g$ is a Cartan subalgebra, and $\Delta$ is a root system.
 \item Simple roots are $\Pi$.
 \item $\Delta^+$ are positive roots with respect to $\Pi$.
 \item for $\alpha \in \Delta^+$, $h_{\alpha}\in \h$ is the corresponding coroot and $e_{\alpha}, f_{\alpha}$ are bases of root subspaces~$\g_{\alpha}$,~$\g_{-\alpha}$
\end{itemize}

\section{Preliminaries}\label{preliminaries}

\subsection[Affine vertex algebra V\^{}kappa(sl\_2) and its modules]{Affine vertex algebra $\boldsymbol{V^{\kappa}(\sl_2)}$ and its modules}

Following \cite{FB,K2,LL}, let us first recall some background on affine Lie algebras and vertex algebras.
Let ${\g}$ be a finite-dimensional simple Lie algebra over $\mathbb{C}$, and let $(\cdot,\cdot)$ be a nondegenerate symmetric bilinear form on ${\g}$.
We fix a triangular decomposition
\[
 {\g} = {\n}_{-} \oplus {\h} \oplus {\n}_{+}.
\]

The affine Lie algebra $\widehat{\g}$ associated with ${\g}$ is defined by
\[
 \widetilde{\g} = {\g} \otimes \mathbb{C}\bigl[t,t^{-1}\bigr] \oplus \mathbb{C}K \oplus \mathbb{C}d,
\]
where $K$ is the canonical central element \cite{K}.
The Lie bracket on $\widetilde{\g}$ is given by
\begin{gather*}
 [x \otimes t^n, y \otimes t^m] = [x,y] \otimes t^{n+m} + n (x,y) \delta_{n+m,0} K, \\ %\label{comut-af-1}\\
 [d, x \otimes t^n] = n x \otimes t^n,  %\label{comut-af-2}
\end{gather*}
for $x,y \in {\g}$ and $m,n \in \mathbb{Z}$.
We write $x(n)$ for $x \otimes t^n$.

The Cartan subalgebra ${\widetilde{\g}}_0$ and subalgebras ${\widetilde{\g}}_{\pm}$ of ${\widehat{\g}}$ are defined by
\[
 {\widetilde{\g}}_0 = {\h} \oplus \mathbb{C}K \oplus \mathbb{C}d, \qquad
 {\widetilde{\g}}_{\pm} = \g \otimes t^{\pm}\mathbb{C}[t^{\pm}] + {\n}_{\pm} \otimes \mathbb{C}.
\]
Thus, we obtain the triangular decomposition
\[
 \widetilde{\g} = \widetilde{\g}_{+} \oplus \widetilde{\g}_{0} \oplus \widetilde{\g}_{-}.
\]

Let $\hg = [ \, \widetilde{\g}, \widetilde{\g} \, ] = \hg_{+} \oplus \hg_{0} \oplus \hg_{-}$,
where $\hg_{0} = {\h} \oplus \mathbb{C}K$ and $\hg_{\pm} = \widetilde{\g}_{\pm}$.
Let
\[
 P = {\g} \otimes \mathbb{C}[t] \oplus \mathbb{C}K.
\]
For each $\kappa \in \mathbb{C}$, let $\mathbb{C}v_{\kappa}$ be the one-dimensional $P$-module on which
${\g} \otimes \mathbb{C}[t]$ acts trivially and~$K$ acts as multiplication by $\kappa$.
Define the generalized Verma module
\[
 V^{\kappa}(\g) = U(\widehat{\g}) \otimes_{U(P)} \mathbb{C}v_{\kappa}.
\]
Then $V^{\kappa}(\g)$ carries a natural structure of a vertex algebra generated by the fields
\[
 x(z) = Y(x(-1)\mathbf{1}, z) = \sum_{n \in \mathbb{Z}} x(n) z^{-n-1},
 \qquad (x \in \g),
\]
where $\mathbf{1} = 1 \otimes v_{\kappa}$ is the vacuum vector.

A $\hg$-module $N$ is called \emph{restricted} (or \emph{smooth}) if for every $w \in N$ and $x \in \g$ we have
\[
 x(z) w \in \mathbb{C}((z)).
\]
A restricted $\hg$-module of level $\kappa$ naturally acquires the structure of a $V^{\kappa}(\g)$-module.

From now on, we specialize to $\g = \mathfrak{sl}_2$ with standard generators $e$, $f$, $h$,
and let ${\frak b} = \mathbb{C}h + \mathbb{C}e$.
Recall that
\[
 [e,f] = h, \qquad [h,e] = 2e, \qquad [h,f] = -2f, \qquad (e,f) = 1, \qquad (h,h) = 2.
\]
Hence $\h = \mathbb{C}h$, $\n_{+} = \mathbb{C}e$, and $\n_{-} = \mathbb{C}f$.
Let $\hb$ denote the subalgebra of $\hg$ generated by $e(n)$ and $h(n)$ for $n \in \mathbb{Z}$.

Assume first that $\kappa \neq -2$.
Define the canonical Sugawara Virasoro vector in the vertex algebra $V^{\kappa}(\mathfrak{sl}_2)$ by
\[
 \omega = \frac{1}{2(\kappa + 2)} \bigl( e(-1)f(-1) + f(-1)e(-1) + \tfrac{1}{2}h(-1)^2 \bigr)\mathbf{1}.
\]
Then the corresponding field
\[
 Y(\omega,z) = L(z) = \sum_{n \in \mathbb{Z}} L(n) z^{-n-2}
\]
satisfies the commutation relations of the Virasoro algebra with central charge
\[
 c_{\kappa} = \frac{3\kappa}{\kappa + 2}.
\]
Thus, every $V^{\kappa}(\mathfrak{sl}_2)$-module becomes a module over the Virasoro algebra.
Moreover,
\begin{equation*} %\label{aff-vir}
 [L(n), x(m)] = -m x(n+m), \qquad x \in \{e,f,h\}.
\end{equation*}
In particular,
\[
 [L(n), x(0)] = 0, \qquad x \in \{e,f,h\}.
\]

Now let $\kappa = -2$.
Set
\[
 t = \tfrac{1}{2}\bigl( e(-1)f(-1) + f(-1)e(-1) + \tfrac{1}{2}h(-1)^2 \bigr)\mathbf{1}
 \in V^{-2}(\mathfrak{sl}_2),
\]
and define
\[
 T(z) = Y(t,z) = \sum_{n \in \mathbb{Z}} T(n) z^{-n-2}.
\]
Then
\[
 [T(n), x(m)] = 0 \qquad \text{for all } m,n \in \mathbb{Z},
\]
so the modes $T(n)$ are central elements.
In particular, $t$ generates the center of the vertex algebra $V^{-2}(\mathfrak{sl}_2)$
(cf.~\cite{efren,FB}).
This center is a commutative vertex algebra $M_T(0)$, which, as a~vector space, is isomorphic to the polynomial algebra
$ \mathbb{C}[ T(-n) \mid n \ge 0 ]$.

\subsection{Weyl vertex algebra and its Whittaker modules}
\label{Weyl}

We now recall the definition of the Weyl vertex algebra, which will serve as one of our main objects of study.
To define this vertex algebra, we first review the notion of the Weyl algebra~$\widehat{\mathcal{A}}$.

Let $\mathcal{L}$ be the infinite-dimensional Lie algebra with generators
$K$, $a(n)$, $a^{*}(n)$, $n \in \mathbb{Z}$,
such that $K$ is central and the only nontrivial commutation relations are
\[
 [a(n), a^{*}(m)] = \delta_{n+m,0} K, \qquad n,m \in \mathbb{Z}.
\]
The corresponding Weyl algebra $\widehat{\mathcal{A}}$ is defined by{\samepage
\[
 \widehat{\mathcal{A}} = \frac{U(\mathcal{L})}{\langle K - 1 \rangle},
\]
where $\langle K - 1 \rangle$ denotes the two-sided ideal generated by $K - 1$.
Hence, in $\widehat{\mathcal{A}}$ we have $K = 1$.}

To construct the Weyl vertex algebra, we consider the simple $\widehat{\mathcal{A}}$-module $W$ generated by a~cyclic vector $\mathbf{1}$ satisfying
\[
 a(n)\mathbf{1} = a^{*}(n+1)\mathbf{1} = 0, \qquad n \ge 0.
\]
Thus,
\[
 W= \mathbb{C}[ a(-n),\, a^{*}(-m) \mid n > 0,\, m \ge 0 ].
\]

By the generating fields theorem, there exists a unique vertex algebra $(W, Y, \mathbf{1})$
such that the vertex operator map
\[
 Y\colon \ W\longrightarrow \operatorname{End}(W)\bigl[\bigl[z, z^{-1}\bigr]\bigr]
\]
is determined by
\[
 Y(a(-1)\mathbf{1}, z) = a(z), \qquad Y(a^{*}(0)\mathbf{1}, z) = a^{*}(z),
\]
where
\[
 a(z) = \sum_{n \in \mathbb{Z}} a(n) z^{-n-1},
 \qquad
 a^{*}(z) = \sum_{n \in \mathbb{Z}} a^{*}(n) z^{-n}.
\]

We choose the following conformal vector (cf.~\cite{KR}) of central charge $c = 2$:
\[
 \omega = \frac{1}{2}\bigl(a(-1)a^{*}(-1) \bigr)\mathbf{1}.
\]
Then $(W, Y, \mathbf{1}, \omega)$ has the structure of a
$ \mathbb{Z}_{\ge 0}$-graded vertex operator algebra.

 Following \cite{ALPY},
we define the Whittaker module for $\widehat{\mathcal{A}}$ to be the quotient \[M_1(\bm{\lambda}, \bm{\mu}) = \widehat{\mathcal{A}}/I,\]
 where $\bm{\lambda} = (\lambda_n)_{n \in {\Z}_{\ge 0}} $ and $\bm{\mu} = (\mu_n) _{n \in {\Z}_{>0}} $ are sequences such that $ \lambda_n=0$ and $\mu_n=0$ for $n \gg 0$ and $I$ is the left ideal
\[ I = \langle \{ a(i) - \lambda_i,\, a^{*}(j)-\mu_j \mid i,j \in {\Z}_{\ge 0},\, j >0\} \rangle.\]

 Let $\mathfrak n$ be the subalgebra of $\mathcal L$ generated by $a(n)$, $a^*(n+1)$, $n \in {\Z}_{\ge 0}$.
 Then $\mathfrak n$ is a commutative, and therefore nilpotent subalgebra of $\mathcal L$.

\begin{Proposition}[\cite{ALPY}] \label{standard} We have
\begin{enumerate}\itemsep=0pt
\item[$1.$] $M_1( \bm{ \lambda}, \bm{ \mu}) $ is a universal Whittaker module for the Whittaker pair $(\mathcal L, \mathfrak n)$. It is of level $K=1$ with the
Whittaker function $\Lambda = (\bm{\lambda}, \bm{\mu}) \colon \mathfrak n \rightarrow {\C}$:
\[
\Lambda(a(i)) = \lambda_i, \qquad \Lambda(a^*(j)) = \mu_j, \qquad i,j \in {\Z}_{\ge 0},\ j>0.
\]
\item[$2.$] $M_1( \bm{ \lambda}, \bm{ \mu}) $ is an irreducible $\widehat{\mathcal A}$-module.
\item[$3.$] $M_1( \bm{ \lambda}, \bm{ \mu}) $ is an irreducible weak module for the Weyl vertex algebra~$W$.
\end{enumerate}
\end{Proposition}

We denote the associated Whittaker vector by ${\bf w}_{ \bm{ \lambda}, \bm{ \mu} }$.

\section[Futorny--Guo--Xue--Zhao modules \protect{[17]}]{Futorny--Guo--Xue--Zhao modules \cite{FGXZ}}
\label{FGXZ-approach}

Let us consider the case $\mathfrak{g} = \mathfrak{sl}_{2}$.
Take integers $N, M \in \mathbb{Z}$ such that $N + M \ge 0$.
In~\cite{FGXZ}, the authors studied the Lie subalgebra
$\mathcal{S}_{N,M} \subset \widehat{\mathfrak{g}}$ generated by $K$ and the elements
\[
 \{ e(n),\, f(p),\, h(m) \mid n > N,\; p > M,\; m \ge 0 \}.
\]
Let $\varphi \colon \mathcal{S}_{N,M} \to \mathbb{C}$ be a Lie algebra homomorphism such that $\kappa = \varphi(K)$.
Denote by $\mathbb{C}v$ the one-dimensional $\mathcal{S}_{N,M}$-module on which $x \in \mathcal{S}_{N,M}$ acts as
\[
 x \cdot v = \varphi(x) v.
\]

We then form the induced $\widehat{\mathfrak{g}}$-module
\[
 \widehat{M}(\varphi) = \operatorname{Ind}_{\mathcal{S}_{N,M}}^{\widehat{\mathfrak{g}}} \mathbb{C}v.
\]
The module \smash{$\widehat{M}(\varphi)$} is a restricted (i.e., smooth) $\widehat{\mathfrak{g}}$-module of level $\kappa$.
Hence, \smash{$\widehat{M}(\varphi)$} is a weak $V^{\kappa}(\mathfrak{g})$-module.
Let $u = 1 \otimes v$.

The following results were proved in~\cite[Theorem~1.1]{FGXZ}.

\begin{Theorem}[\cite{FGXZ}]\label{fgxz}
Let $\mathfrak{g} = \mathfrak{sl}_{2}$ and $\kappa \ne -2$.
Then \smash{$\widehat{M}(\varphi)$} is an irreducible $V^{\kappa}(\mathfrak{g})$-module
if and only if $\varphi(h(N + M + 1)) \ne 0$.
\end{Theorem}

In~\cite[Theorem~1.2]{FGXZ}, the authors considered the critical level and, in certain cases,
described the irreducible quotients of \smash{\smash{$\widehat{M}(\varphi)$}} when $\varphi(h(N+M+1)) \ne 0$. We shall slightly reformulate their theorem.

Let
$
 T(z) = \sum_{n \in \mathbb{Z}} T(n) z^{-n-2}
$
denote the center of $V^{-2}(\mathfrak{g})$ at the critical level. Note that
on every irreducible $V^{-2}(\mathfrak{g})$-module, the central elements $T(n)$ must act as scalar.

By construction of the module \smash{$\widehat{M}(\varphi)$}, one can see that for $N_0 = N+M+1 \in {\Z}_{>0}$ and $i\in {\Z}_{> 0}$ we have
\begin{itemize}\itemsep=0pt
\item $T(N_0-i)u, u$ are linearly independent.
\item $T(N_0+ i-1)u$ is proportional to $u$ and $T(N_0+ i-1)u =0$ for $i \gg 0$.
 \end{itemize}
 Choose
\[ \theta(z) = \sum_{n\in {\Z}} \theta_n z^{-n-2} \in {\C} ((z)) \]
 such that $T(N_0+ i-1)u = \theta_{N_0+i-1} u$ for $i > 0$.
 (Other coefficients are arbitrary complex numbers.)

 Recall that a nonzero vector $\widehat{M}(\varphi, \theta)$ is called a singular vector of weight $\lambda\in\C$ if
\[
e(N+i)w=f(M+i)w=(h(i)-\varphi(h(i)))w=0, \qquad i\in\Z_{>0},
\]
and
\[
h(0)w=\lambda w.
\]
Since $T(n)$ is central at the critical level, it commutes with $e(N+j)$, $f(M+j)$, and $h(j)$ for all $j>0$. Hence, by the choice of $\theta(z)$, the vectors
\begin{equation}\label{singul}
(T(N_0-i)-\theta_{N_0-i})u, \qquad i>0,
\end{equation}
are singular vectors.

Let $\widehat{M}(\varphi, \theta)$ be the submodule of \smash{$\widehat{M}(\varphi)$} generated by the singular vectors (\ref{singul}). Define the quotient module
 \[ \overline M(\varphi, \theta) = \widehat M(\varphi)/ \widehat{M}(\varphi, \theta). \]

The following result is a slightly reformulated version of \cite[Theorem 1.2]{FGXZ}.

\begin{Theorem}[\cite{FGXZ}]\label{fgxz-critical}
Let $\mathfrak{g} = \mathfrak{sl}_{2}$ and $\kappa = -2$.
Suppose that $N+ M \in {\Z}_{\ge 0}$, $\varphi(h(N + M + 1)) \ne 0$ and let $\theta(z)$ is as above. Then
 $\overline{M}(\varphi, \theta)
$
is an irreducible $V^{\kappa}(\mathfrak{g})$-module such that
 \[ T(z) \equiv \theta(z) \qquad \mbox{on} \quad \overline{M}(\varphi, \theta).\]
\end{Theorem}

We shall see below that these results are closely related to the irreducible Wakimoto modules
at the critical level constructed in~\cite{A-07, A-14}.
The results of~\cite{FGXZ} do not describe the irreducible quotients of the universal modules
\smash{$\widehat{M}(\varphi)$} in the cases where these modules are reducible.
However, the paper~\cite{A-14} already provided a complete classification of all irreducible quotients
of \smash{$\widehat{M}(\varphi)$} at the critical level.
In the present paper, we also aim to describe the simple quotients for non-critical levels.

Consider the following automorphism of $\widehat{\mathfrak{sl}}_{2}$:
\[
 e(n) \mapsto f(n), \qquad
 f(n) \mapsto e(n), \qquad
 h(n) \mapsto -h(n), \qquad
 K \mapsto K.
\]
This automorphism extends to an automorphism $\tau$ of order two on $V^{\kappa}(\mathfrak{sl}_{2})$.
Note that the subalgebra $\mathcal{S}_{N,N}$ of $\widehat{\mathfrak{sl}}_{2}$ is $\tau$-invariant.
Let $V^{\kappa}(\mathfrak{sl}_{2})^{\langle \tau \rangle}$ denote the fixed-point subalgebra under~$\tau$.

\begin{Lemma}
Let $\varphi \colon \mathcal{S}_{N,N} \to \mathbb{C}$ be a Lie functional such that
$\varphi(e(n)) = \varphi(f(n)) = 0$ for $n > N$,
and $\varphi(h(m)) = 0$ for sufficiently large~$m$, but $\varphi \ne 0$.
Then
\[
 \tau \circ \widehat{M}(\varphi) = \widehat{M}(-\varphi).
\]
In particular, $\widehat{M}(\varphi) \not\cong \tau \circ \widehat{M}(\varphi)$.
\end{Lemma}

Using the main result of~\cite{ALPY}, we obtain the following.

\begin{Corollary}
Assume that $\kappa \ne -2$. The module \smash{$\widehat{M}(\varphi)$} is an irreducible
$V^{\kappa}(\mathfrak{sl}_{2})^{\langle \tau \rangle}$-module.
\end{Corollary}

\section[Wakimoto modules for widehat sl\_2]{Wakimoto modules for $\boldsymbol{\widehat{\mathfrak{sl}}_{2}}$}
\label{section-Wakimoto-1}

Let $\mathfrak{h} = \mathbb{C}b$ be a one-dimensional commutative Lie algebra equipped with a symmetric bilinear form
$(b,b) = 2$, and let $\widehat{\mathfrak{h}} = \mathfrak{h} \otimes \mathbb{C}\bigl[t,t^{-1}\bigr] \oplus \mathbb{C}c$
be its affinization.
Set $b(n) = b \otimes t^{n}$.
Let $\pi^{\kappa+2}$ denote the simple $\widehat{\mathfrak{h}}$-module of level $\kappa+2$ generated by the vector ${\bf 1}$ such that
\[
 b(n) {\bf 1} = 0 \qquad \text{for all } n \ge 0.
\]
As a vector space, we have
\[
 \pi^{\kappa+2} = \mathbb{C}[b(n) \mid n \le -1].
\]
Then $\pi^{\kappa+2}$ carries a unique structure of a vertex algebra generated by the field
\[
 b(z) = \sum_{n \in \mathbb{Z}} b(n) z^{-n-1},
\]
satisfying the commutation relation
\[
 [b(n), b(m)] = 2(\kappa + 2) n \delta_{n+m,0}.
\]

Let $V^{\kappa}(\mathfrak{sl}_{2})$ denote the universal vertex algebra of level~$\kappa$ associated to the affine Lie algebra~$\widehat{\mathfrak{sl}}_{2}$.
Recall that $W$ is the Weyl vertex algebra from Section~\ref{Weyl}.
There exists an injective homomorphism of vertex algebras
\[
 \Phi \colon \ V^{\kappa}(\mathfrak{sl}_{2}) \hookrightarrow W \otimes \pi^{\kappa+2},
\]
given on generators by
\begin{align*}
 e &= a(-1){\bf 1}, \\
 h &= -2 a^{*}(0)a(-1){\bf 1} + b(-1){\bf 1}, \\
 f &= -a^{*}(0)^{2} a(-1){\bf 1} + \kappa a^{*}(-1){\bf 1} + a^{*}(0) b(-1){\bf 1}.
\end{align*}

We identify
\[
 a = a \otimes {\bf 1}, \qquad
 a^{*} = a^{*} \otimes {\bf 1}, \qquad
 b = {\bf 1} \otimes b.
\]
For $x \in \{ e, f, h \}$, set
\[
 x(z) = \sum_{n \in \mathbb{Z}} x(n) z^{-n-1}.
\]
Then
\begin{align*}
 e(z) &= Y(e,z) = a(z), \\
 h(z) &= Y(h,z) = -2 {:}a^{*}(z)a(z){:} + b(z), \\
 f(z) &= Y(f,z) = -{:}a^{*}(z)a^{*}(z)a(z){:} + \kappa \partial_{z}a^{*}(z) + a^{*}(z)b(z).
\end{align*}

The following proposition is a standard result in the theory of vertex algebras
(cf.~\cite{FB, K2, LL}), applied to the vertex algebra $W \otimes \pi^{\kappa+2}$.

\begin{Proposition}\label{constr-new1}
Assume that $M_{1}$ is a restricted module for the Weyl algebra and $N_{1}$ is a~restricted module of level $\kappa+2$ for the Heisenberg algebra $\widehat{\mathfrak{h}}$.
That is, for every $u \in M_{1}$ and $v \in N_{1}$ there exists $N \in \mathbb{Z}_{\ge 0}$ such that
\[
 a(n)u = 0, \qquad b(n)v = 0 \qquad \text{for all } n \ge N.
\]
Then $M_{1} \otimes N_{1}$ is an $W \otimes \pi^{\kappa+2}$-module, and therefore a
$V^{\kappa}(\mathfrak{sl}_{2})$-module.
\end{Proposition}

Assume now that $\kappa \ne -2$.
Then we have a natural action of the Virasoro algebra generated by the Sugawara vector
\begin{align*}
 \omega
 &= \frac{1}{2(\kappa + 2)}
 \big( e(-1)f(-1) + f(-1)e(-1) + \tfrac{1}{2}h(-1)^{2} \big){\bf 1} \\
 &= a(-1)a^{*}(-1)
 + \frac{1}{4(\kappa + 2)} \big( b(-1)^{2} - 2b(-2) \big).
\end{align*}
Let
\[
 L(z) = Y(\omega, z) = \sum_{n \in \mathbb{Z}} L(n) z^{-n-2}.
\]
For further details, see~\cite{FB}.

\section[Wakimoto realization of modules widehat M(varphi): The case kappa not = -2]{Wakimoto realization of modules $\boldsymbol{\widehat M(\varphi)}$: The case $\boldsymbol{\kappa \ne -2}$}

Let $\widehat{\mathfrak h}_{\ge 0} = {\mathfrak h} \otimes {\C}[t]$. Let $\bm{\eta} \in (\widehat{\mathfrak h}_{\ge 0}) ^* $ such that
 \[ \bm{\eta}(b (n)) = 0, \qquad n \gg 0. \]
Let $N_1 (\bm{\eta}) $ be the standard Whittaker module for the Heisenberg Lie algebra $\widehat{\mathfrak h}$ of level $\kappa +2$ with the Whittaker function $\bm{\eta}$. This module is an irreducible module for the Heisenberg vertex algebra $\pi^{\kappa+2}$ (for details see \cite{ALPY} and references therein).

Now we shall prove the first main theorem in our paper.

\begin{Theorem} \label{prvi-teorem} Assume that $N >0$ such that $\bm{\eta}(b(N)) \ne 0$ and $\bm{\eta}(b(j)) =0$ for $j >N$.
Then $W \otimes N_1(\bm{\eta})$ is an irreducible $V^{\kappa}(\sl_2)$-module. Moreover,
 \[ W \otimes N_1(\bm{\eta}) \cong \widehat{M}(\varphi), \]
where $\varphi \colon S_{-1, N} \rightarrow \C$ is the Lie functional uniquely determined by
 \[ \varphi(e(n)) = \varphi(f(N+1+n)) = 0, \qquad \varphi(h(n)) = \bm{\eta}(b(n)), \qquad n \in {\Z}_{\ge 0}. \]
\end{Theorem}
\begin{proof}
Identify $v_{\bm{\eta}}$ with $ {\bf 1} \otimes v_{\bm{\eta}}$.
We first claim that
\begin{itemize}\itemsep=0pt
\item[(1)] $v_{\bm{\eta}}$ is a cyclic vector for the $V^{\kappa}(\sl_2)$-action, i.e.,
 \[ V^{\kappa}(\sl_2). v_{\bm{\eta}} = W \otimes N_1(\bm{\eta}). \]
\end{itemize}
Let us prove claim (1). The PBW basis of $W \otimes N_1(\bm{\eta})$ as $W \otimes \pi^{\kappa +2}$-module is
\[ a(-n_1) \cdots a(-n_r) a^*(-m_1) \cdots a^* (-m_s) b(-j_1) \cdots b(-j_t) v_{\bm{\eta}}, \]
where $r,s,t \in {\Z}_{\ge 0}$, $n_1, \dots, n_r, j_1,\dots, j_t \in {\Z}_{>0}$, $m_1, \dots, m_s \in {\Z}_{\ge 0}$.
Define
\[
\begin{aligned}
&\deg\bigl(a(-n_1)\cdots a(-n_r)a^*(-m_1)\cdots a^*(-m_s)
 b(-j_1)\cdots b(-j_t)\bigr)\\
&\qquad = n_1+\cdots+n_r+m_1+\cdots+m_s+j_1+\cdots+j_t.
\end{aligned}
\]
Then
\begin{gather*}
 W \otimes N_1(\bm{\eta}) = \bigoplus_{m \in {\Z}_{\ge 0} } ( W \otimes N_1(\bm{\eta}) )_m, \\
w \in ( W \otimes N_1(\bm{\eta}) )_m \iff \deg (w) = m.
\end{gather*}
We shall prove by induction on $m$ that $( W \otimes N_1(\bm{\eta}) )_m \subset V^{\kappa}(\sl_2). v_{\bm{\eta}}$.

Basis: $m = 0$.
Then each $w \in ( W \otimes N_1(\bm{\eta}) )_0$ is a linear combination of vectors
$a^*(0) ^p v_{\bm{\eta}}$.

Using the action of $f(N)$ on $W \otimes N_1(\bm{\eta})$, we easily get \[ f(N) ^p v_{\bm{\eta}} = \nu a^*(0) ^p v_{\bm{\eta}}, \qquad (\nu \ne 0), \]
which implies that $a^*(0) ^p v_{\bm{\eta}} \in V^{\kappa}(\sl_2). v_{\bm{\eta}}$ for each $p \in {\Z}_{\ge 0}$. Therefore,
$( W \otimes N_1(\bm{\eta}) )_0 \subset V^{\kappa}(\sl_2). v_{\bm{\eta}}$. The basis of induction holds.

Assume now that there is $M \in {\Z}_{>0}$ such that
 \[( W \otimes N_1(\bm{\eta}) )_m \subset V^{\kappa}(\sl_2). v_{\bm{\eta}} \qquad \mbox{for} \ m < M. \]
Let $w \in ( W\otimes N_1(\bm{\eta}) )_M$. It suffices to consider the case when $w $ is written in PBW basis. So assume that
\begin{align}
&w=a(-n_1)\cdots a(-n_r)a^*(-m_1)\cdots a^*(-m_s)b(-j_1)\cdots b(-j_t)v_{\bm{\eta}}, \label{basis-11} \\
&M= n_1 + \cdots + n_r + m_1 +\cdots + m_s + j_1 +\cdots j_t. \nonumber
\end{align}
If $r >0$, then we have $ w= e(-n_1) w'$ for $\deg(w') < M$. By the inductive assumption, we have that $w ' \in V^{\kappa}(\sl_2). v_{\bm{\eta}}$, which implies that $w \in V^{\kappa}(\sl_2). v_{\bm{\eta}}$.
So it suffices to prove the statement in the case $r=0$. Therefore, let
 \[w= a^*(-m_1) \cdots a^* (-m_s) b(-j_1) \cdots b(-j_t) v_{\bm{\eta}}, \qquad
  M= m_1 +\cdots + m_s + j_1 +\cdots +j_t. \]
Assume that $s >0$. Then
 \[ f(N-m_1) a^*(-m_2) \cdots a^* (-m_s) b(-j_1) \cdots b(-j_t) v_{\bm{\eta}} = \nu_1 w + w', \]
for certain $\nu_1 \ne 0$ and $\deg(w') <M$. Using the inductive assumption again, we get that $w \in V^{\kappa}(\sl_2). v_{\bm{\eta}}$.

Finally, we consider the case $r= s=0$, and $w$ of the form
 \begin{align*}
 w= b(-j_1) \cdots b(-j_t) v_{\bm{\eta}}, \qquad
 M= j_1 +\cdots +j_t.
 \end{align*}
Let $t >0$. Then
 \[ h(-j_1) b(-j_2) \cdots b(-j_t) v_{\bm{\eta}} = w + w_1, \]
where $\deg(w_1) = M $ and $w_1$ is a linear combination of basis elements of the form (\ref{basis-11}), where ${r >0}$. We have already shown that in these cases $w_1 \in V^{\kappa}(\sl_2). v_{\bm{\eta}}$. This proves that $w \in V^{\kappa}(\sl_2). v_{\bm{\eta}}$. Therefore, the claim (1) holds.

So $v_{\bm{\eta}}$ is cyclic vector in $ W\otimes N_1({\bm{\eta}})$. Since
 \[ e(n) v_ {\bm{\eta}} = f(N+1 + n) v_ {\bm{\eta}} =0, \qquad h(n) v_ {\bm{\eta}} = {\bm{\eta}} (b(n)) v_ {\bm{\eta}}, \]
 for $n\geq 0$,
we conclude that $ W\otimes N_1({\bm{\eta}})$ is isomorphic to a quotient of the universal module \smash{$\widehat{M}(\varphi)$}, where
$ \varphi\colon \mathcal S_{-1, N} \rightarrow {\C}$ is the Lie functional defined by
 \[ \varphi(e(n)) = \varphi(f(N+1+ n)) =0, \qquad \varphi(h(n)) = {\bm{\eta}} (b(n)), \qquad n \in {\Z}_{\ge 0}. \]
Since $ \varphi(h(N)) \ne 0$, Theorem~\ref{fgxz} implies that \smash{$\widehat{M}(\varphi)$} is irreducible. Therefore,
 \[ \widehat{M}(\varphi) = W \otimes N_1({\bm{\eta}}), \]
and thus $W\otimes N_1({\bm{\eta}})$ is an irreducible $ V^{\kappa}(\sl_2)$-module.
\end{proof}

\section[The explicit realization of simple quotients of widehat M(varphi) at the critical level]{The explicit realization of simple quotients of $\boldsymbol{\widehat M(\varphi)}$\\ at the critical level}
\label{ired-critical-sl(2)}

Now we consider the critical level case $\kappa =-2$. Therefore, $\pi^{\kappa +2} = \pi^0 $ is a commutative vertex algebra generated by the field
 \[ b(z) = \sum_{n \in {\Z}} b(n) z^{-n-1}, \qquad b(n)=0 \qquad \mbox{for} \ n \ge 0. \]

Let \[ \bm{\chi}(z) = \sum_{n \in {\Z} } \chi_n z^{-n-1} \in {\C}((z)). \]
Let $L (\chi) $ be the one-dimensional $\pi^0$-module such that $b(n)$ acts as a multiplication with $\chi_n$. Since $\pi^0$ is commutative, all irreducible modules are one-dimensional and isomorphic to $L (\bm{\chi}) $ for certain $\bm{\chi}(z) \in {\C}((z))$.

Note that $L (\bm{\chi})$ is a simple quotient of the Whittaker $\pi^0$-module $N_1(\bm{\eta})$ where
$\bm{\eta} \in (\widehat{\mathfrak h}_{\ge 0}) ^* $ is given by
 \[ \bm{\eta}(b (i)) = \chi_i, \qquad i \ge 0. \]

Now we consider the Wakimoto module
$W_{\bm{\chi}} = W \otimes L (\bm{\chi})$.
The irreducibility result for modules $W_{\chi}$ was obtained in \cite{A-07, A-14}.
In order to present the irreducibility criterion, we need to
recall the definition of Schur polynomials.

 Define the Schur polynomials $S_{r}(x_{1},x_{2},\dots)$
 in variables $x_{1},x_{2},\dots$ by the following equation:
\begin{equation*}%\label{eschurd}
\exp\left(\sum_{n=1}^{\infty}\frac{x_n}{n}y^n\right)
=\sum_{r=0}^{\infty}S_r(x_1,x_2,\dots)y^r.
\end{equation*}

We shall also use the following formula for Schur polynomials:
\begin{equation*}%\label{det-schur}
S_r(x_1,x_2,\dots)=\frac{1}{r!}
\det\begin{pmatrix}
 x_1 & x_2 & \cdots & & x_r \\
 -r+1 & x_1 & x_2 & \cdots & x_{r-1} \\
 0 & -r+2 & x_1 & \cdots & x_{r-2} \\
 0 & \ddots & \ddots & \ddots & \\
 0 & \cdots & 0 & -1 & x_1
\end{pmatrix}.
\end{equation*}

\begin{Theorem}[\cite{A-07, A-14}] \label{A14}
Assume that $\bm{\chi} \in {\C}((z))$. Then the Wakimoto module
$W_{-\bm{\chi}}$ is an irreducible $V^{-2}(\mathfrak{sl}_2)$-module if and only
if $\bm{\chi}$ satisfies one of the following conditions:
\begin{enumerate}\itemsep=0pt
\item[$(i)$] There is $p \in {\Z}_{>0}$
 such that
\[ \bm{\chi}(z) = \sum_{n=-p} ^{\infty} {\chi}_{-n} z ^{n-1} \in
{\C}((z)) \qquad \mbox{and} \qquad
\chi_p \ne 0. \]

\item[$(ii)$] \[ \bm{\chi}(z) = \sum_{n=0} ^{\infty} {\chi}_{-n} z
^{n-1} \in {\C}((z)) \qquad \mbox{and} \qquad
\chi_0 \in \{1 \} \cup ({\C} \setminus {\Z}). \]
\item[$(iii)$] There is $\ell \in {\N}$ such that
\[ \bm{\chi}(z) = \frac{\ell +1}{z} + \sum_{n=1} ^{\infty} {\chi}_{-n} z ^{n-1} \in
{\C}((z))\]
and $S_{\ell}(-\chi) \ne 0$,
where $S_{\ell}(-\chi) = S_{\ell}(-\chi_{-1}, -\chi_{-2}, \dots)$
is a Schur polynomial.
\end{enumerate}
\end{Theorem}

Let us now explicitly compare modules $\overline M(\varphi, \theta )$ from \cite{FGXZ}, which we recalled in Section~\ref{FGXZ-approach}, with the Wakimoto modules considered in \cite{A-07, A-14}.

\begin{Theorem} \label{krit-kvocijenti}
 Assume that $N=p \in {\Zp}$ and
 \[ \bm{\chi}(z) = \sum_{n=-p} ^{\infty} {\chi}_{-n} z ^{n-1} \in
{\C}((z)). \]
Let $\varphi \colon S_{-1, N} \rightarrow {\C}$ be the Lie functional such that \[\varphi(h(i)) = -\chi_i, \qquad \varphi(e(i)) = \varphi(f(i+N+1)) =0 \qquad \mbox{for} \ i \ge 0. \]

Assume that $\bm{\chi} \in {\C}((z))$ satisfies the condition $(i)$, $(ii)$ or $(iii)$ of Theorem~{\rm \ref{A14}}. Then $W_{-\chi}$ is a simple quotient of $\widehat M( \varphi)$.

 If $p\ge 1$ and $\chi_p \ne 0$, then
 $\overline M(\varphi, \theta) \cong W_{-\bm{\chi}}$, where
 \[ \theta(z) = \frac{1}{2} \left( \bm{\chi}(z) ^2 + 2 \frac{{\rm d}}{{\rm d} z}\bm{\chi} (z) \right) = \sum_{i\ge 0} \theta_{2p +i} z^{-2 p- i - 2}. \]
 \end{Theorem}

 \begin{proof}
By construction, $W_{-\bm{\chi}}$ is generated by a cyclic vector $w$ satisfying conditions
 \[ e(i) w = f(i+ N+1) w =0, \ h(i) w = -\chi _i w, \qquad i\ge 0. \]
Therefore, it is isomorphic to a quotient of $\widehat M( \varphi)$.

Moreover, whenever one of the conditions \textnormal{(i)}, \textnormal{(ii)}, or \textnormal{(iii)} of Theorem~\ref{A14} is satisfied, the module $W_{-\boldsymbol{\chi}}$ is simple. Consequently, in this case, $W_{-\boldsymbol{\chi}}$ is a simple quotient of \smash{$\widehat{M}(\varphi)$}.

Assume next that the condition (i) holds, i.e., $p \ge 1$ and $\chi_p \ne 0$. Since we know that $\overline M(\varphi, \theta)$ is a quotient of \smash{$\widehat M( \varphi)$} by a submodule generated by singular vectors
 \[ (T(i) - \theta_i) u, \qquad i < N_0 = N, \]
and $T(z) = \theta(z)$ on $W_{-\chi}$, we get that $W_{-\chi}$ is a quotient of $\overline M(\varphi, \theta)$.

For $N\ge 1$, the module $\overline M(\varphi, \theta)$ is irreducible by Theorem~\ref{fgxz-critical}, so $\overline M(\varphi, \theta) \cong W_{-\bm{\chi}}$.
 \end{proof}

Let us now consider the case $N=0$. The simple quotients in the case $N=0$ were not described in~\cite{FGXZ}. So our result in the case $N=0$ is a refinement of the main result of~\cite{FGXZ}.
Moreover, we shall see below that in this case $\overline M(\varphi, \theta)$ need not be irreducible.
 Applying previous theorem and \cite[Corollary 5.2 and Remark 5.1]{A-14}, we get the following.

\begin{Proposition} \label{reducibilni}
 Assume that \[ \bm{\chi}(z) = \frac{\ell +1}{z} + \sum_{n=1} ^{\infty} {\chi}_{-n} z ^{n-1} \in
{\C}((z))\]
is such that
 $\ell \in {\Z}_{>1} $ and $S_{\ell}(-\chi) =0 $ or $\ell <0$.
Then $\overline M(\varphi, \theta)$ is reducible and contains an irreducible submodule isomorphic to $W_{-\chi} ^{{\rm int}_{\mathfrak{sl}_2}}$, which is the maximal integrable $\mathfrak{sl}_2$-submodule of~$W_{-\chi}$.
 \end{Proposition}

\begin{Remark}
The reducibility was established in the earlier paper \cite{A-14}, using methods developed in \cite{A-07}. These methods are based on the vertex superalgebra $\mathcal V$, which appears as a limit of the $N=2$ superconformal vertex algebra, together with the critical-level version of the correspondence between $N=2$ superconformal vertex algebra and affine vertex algebras associated to $\mathfrak{sl}_2$, and in particular on the anti-Kazama--Suzuki mapping (cf. \cite{A-1999, FST}). Since it is much easier to prove irreducibility or reducibility for certain $\mathcal V$-modules, one can then deduce the corresponding irreducibility or reducibility of $\widehat{\mathfrak{sl}}_2$-modules.

More precisely, in \cite{A-07} we constructed $\mathcal V$-modules $\widetilde F_{\bm{\chi}}$, realized on certain fermionic Fock spaces and parametrized by $\bm{\chi} \in \C((z))$. We also constructed a functor $\mathcal L_0$ from the category of $\mathcal V$\nobreakdash-mod\-ules to the category of $V^{-2}(\mathfrak{sl}_2)$-modules, and showed that
\smash{$W_{-\chi} \cong \mathcal L_0\bigl(\widetilde F_{\bm{\chi}}\bigr)$}.
Consequently, the question of irreducibility or reducibility of Wakimoto modules reduces to the corresponding question for the module \smash{$\widetilde F_{\bm{\chi}}$}. If $\bm{\chi}$ satisfies the condition of Proposition~\ref{reducibilni}, we proved that \smash{$\widetilde F_{\bm{\chi}}$} is reducible and contains an irreducible submodule $U_{\bm{\chi}}$ such that
\smash{$\mathcal L_0(U_{\bm{\chi}}) = W_{-\chi}^{{\rm int}_{\mathfrak{sl}_2}}$}.
In this way, we obtained the reducibility of $W_{-\chi}$.
\end{Remark}

\section{Wakimoto realization of generalized Whittaker modules}
 \label{gen-whitt-sect}

Let $\mathcal P_{q, r}$, where $q, r \in {\Z}_{>0}$, $q < r$, be the Lie subalgebra of $\widehat{\g}$ generated by $K$ and the elements
\[
 \{ e(n),\, h(m),\, f(p) \mid n \ge 0 ,\, m \ge q, \, p \ge r \}.
\]
Let $\psi \colon \mathcal{P}_{q, r} \to \mathbb{C}$ be a homomorphism of the Lie algebra such that $\kappa = \psi(K)$.

Denote by $\mathbb{C}v$ the one-dimensional $\mathcal{P}_{q,r }$-module on which $x \in \mathcal{P}_{q, r}$ acts as
\begin{equation}\label{whit-gen}
x\cdot v=\psi(x) v.
\end{equation}
We then form the induced $\widehat{\mathfrak{g}}$-module
\[
 \widehat{R}(\psi) = \operatorname{Ind}_{\mathcal{P}_{q, r}}^{\widehat{\mathfrak{g}}} \mathbb{C}v.
\]
The module $\widehat{R}(\psi)$ is a restricted (i.e., smooth) $\widehat{\mathfrak{g}}$-module of level $\kappa$.
Hence, $\widehat{R}(\psi)$ is a weak $V^{\kappa}(\mathfrak{g})$-module. This module is a standard (universal) generalized Whittaker module with the Whittaker function $\psi$.

\begin{Definition}Assume that $Z$ is any $V^{\kappa}(\mathfrak{sl}_2)$-module and $v \in Z$. We say that $v$ is a $\psi$\nobreakdash-Whit\-tak\-er vector if it satisfies condition (\ref{whit-gen}).
 \end{Definition}

We shall consider $W \otimes \pi^{\kappa +2}$-modules $M_{1} (\bm{\lambda}, \bm{\mu}) \otimes N_1(\bm{\eta})$ from Section~\ref{Weyl}. Assume that
 \[ \lambda_N \ne 0, \mu_M \ne 0 \qquad \mbox{and} \qquad \lambda_{N+ i} = \mu_{M+i} = 0 \qquad \mbox{for} \ i \in {\Z}_{>0}. \]
As before for $\bm{\eta} \in \bigl(\, \widehat{\mathfrak h}_{\ge 0}\bigr) ^*$, set $\eta_n = \bm{\eta} (b(n))$. Assume that
 \[ \eta_P \ne 0, \qquad \eta_{P+i } =0 \qquad \mbox{for} \ i >0. \]
Using the Frenkel--Feigin homomorphism, we can consider $M_{1} (\bm{\lambda}, \bm{\mu}) \otimes N_1(\bm{\eta})$ as a $V^{\kappa} (\mathfrak{sl}_2)$-module.
Set $w = v_{\bm{\lambda}, \bm{\mu}} \otimes v_{\bm{\eta}}$.
Let
 \[ q = \max \{ M, N+ 1\}, \qquad r = \max \{M+N+1, 2M, P+1\}. \]

By a direct calculation, we see that
\begin{equation*}%\label{uvjeti-1}
e(i)w=A_iw, \qquad h(q+i)w=B_iw, \qquad f(r+i)w=C_iw \qquad i\ge0,
\end{equation*}
where
\begin{align}
&A_i= \lambda_i,\qquad
B_i = -2\sum_{l=1}^M\mu_l\lambda_{q+i-l}+\eta_{q+i},\nonumber\\
&C_i= -\sum_{\substack{k_1,k_2\in {\Z}_{\ge1}\\ k_3\in {\Z}_{\ge0}\\ k_1+k_2+k_3=r+i}}
\mu_{k_1}\mu_{k_2}\lambda_{k_3}
+\sum_{\substack{j_1\in {\Z}_{\ge1}\\ j_2\in {\Z}_{\ge0}\\ j_1+j_2=r+i}}
\mu_{j_1}\eta_{j_2}-\kappa(r+i)\mu_{r+i}.\label{parametri}
\end{align}

 Therefore, $w$ is a $\psi$-Whittaker vector with Whittaker function $\psi \colon \mathcal P_{q,r} \rightarrow \C$ defined by
\begin{equation}\label{uvjeti-2}
\psi(e(i))=A_i, \qquad \psi(h(q+i))=B_i, \qquad \psi(f(r+i))=C_i, \qquad i\ge0,
\end{equation}
 where $A_i$, $B_i$, $C_i$ are given by (\ref{parametri}).

 In this way, we have proved the following.

 \begin{Proposition}
 Let $
\kappa \in {\C}$ be arbitrary. Then $U\bigl(\widehat{\frak{sl}_2}\bigr)w$ is a generalized Whittaker module with the Whittaker function $\psi$ satisfying the condition \eqref{uvjeti-2}.
 \end{Proposition}

 \begin{Remark}
 Using the Frenkel--Feigin realization, the formulas in (\ref{parametri})
are obtained by evaluating the modes of
$
e(z), h(z),
f(z)
$
on the Whittaker vector $w=v_{\lambda,\mu}\otimes v_\eta$.
The choice
\[
q=\max\{M,N+1\}, \qquad r=\max\{M+N+1,2M,P+1\}
\]
ensures that the relevant positive modes act on $w$ by scalars, so that $w$ is a
$\psi$-Whittaker vector for the subalgebra $P_{q,r}$.
More precisely, the coefficients in (\ref{parametri}) come respectively from the linear
term $e(z)$, the quadratic term $-2{:}a^*a{:}$ together with $b(z)$, and the three
summands in the formula for $f(z)$.
 \end{Remark}

\begin{Example} To illustrate the choice of Whittaker function, we recall the following example from \cite[Section~9]{ALZ-16}:
\[
\bm{\lambda}=(\lambda,0,0,\dots), \qquad
\bm{\mu}=(\mu,0,0,\dots), \qquad
\bm{\eta}=(\eta_0,\eta_1,0,\dots),
\]
with $\lambda\neq 0$, $\mu\neq 0$, and $\eta_1\neq 0$. Then
$N=0$, $M=1$, $P=1$,
so
\[
q=\max\{M,N+1\}=1,\qquad r=\max\{M+N+1,2M,P+1\}=2.
\]
The formulas in (\ref{parametri}) reduce to
\begin{gather*}
e(0)w=\lambda w,\qquad e(i)w=0, \qquad i>0,
\\
h(1)w=(\eta_1-2\mu\lambda)w,\qquad h(1+i)w=0, \qquad i>0,
\end{gather*}
and
\[
f(2)w=\mu (\eta_1-\mu \lambda)w,\qquad f(2+i)w=0, \qquad i>0.
\]
Hence $w=v_{\bm{\lambda},\bm{\mu}}\otimes v_{\bm{\eta}}$ is a $\psi$-Whittaker vector for $P_{1,2}$ with
\[
\psi(e(0))=\lambda,\qquad
\psi(h(1))=\eta_1-2\mu\lambda,\qquad
\psi(f(2))=\mu\eta_1-\mu^2\lambda,
\]
and all higher values $\psi(e(i))$, $\psi(h(1+i))$, $\psi(f(2+i))$ vanish for $i>0$.
\end{Example}

 We have the following general conjecture.
 \begin{Conjecture} \label{conj7.3} Assume that $\mu_M\ne 0$, ${\l}_N \ne 0$, $\bm{\eta} \ne 0$ and that $\psi$ is a Whittaker function satisfying \eqref{uvjeti-2}.
 Assume that $\kappa \ne -2$. Then
 $\widehat{R}(\psi)$ is an irreducible $V^{\kappa} (\mathfrak{sl}_2)$-module and
 \[ \widehat{R}(\psi) \cong M_1 (\bm{\lambda}, \bm{\mu}) \otimes N_1 (\bm{\eta}). \]
 \end{Conjecture}

 \begin{Remark}
Although Conjecture~\ref{conj7.3} may be viewed as a natural extension of Theorem~\ref{prvi-teorem}, we do not see a direct way to extend the proof of Theorem~\ref{prvi-teorem} to this setting. The proof of Theorem~\ref{prvi-teorem} uses ingredients specific to Wakimoto-type modules and their free-field realization; in addition, it depends in part on the result of \cite{FGXZ}, for which no direct analogue is presently available for Conjecture~\ref{conj7.3}.

By contrast, the modules arising in Conjecture~\ref{conj7.3} are closer in nature to generalized Whittaker modules, and to the irreducibility arguments used in \cite{ALZ-16}. Thus, while Conjecture~\ref{conj7.3} is a~natural continuation of Theorem~\ref{prvi-teorem} at the level of the expected statement, the relevant methods appear to be different.

This is already reflected at the critical level, where Wakimoto modules and Whittaker-type modules require essentially different proofs. In what follows, we consider the critical-level case, which provides further support for Conjecture~\ref{conj7.3} in view of the analogy with generalized Whittaker modules. \end{Remark}

 Similarly to Theorem~\ref{fgxz-critical} one shows that $ \widehat{R}(\psi)$ is reducible at the critical level. But we shall identify its irreducible subquotients following approaches from \cite{A-14, ALZ-16}.

 Let
 \[ \bm{\chi}(z) = \sum_{n =-P } ^{\infty} \chi_{-n} z^{n-1} \in {\C}((z)), \qquad \chi_P \ne 0, \]
 such that
 \[ \chi_i = \eta_i \qquad \mbox{for} \ i \ge 0, \]
 and let $L(\chi)$ be the simple one-dimensional $\pi^0$-module as in Section~\ref{ired-critical-sl(2)}. Clearly, $L(\bm{\chi})$ is a simple quotient of $N_1 (\bm{\eta})$.
We have the following result.

\begin{Theorem} Assume that $\kappa =-2$. Then $M_1 (\bm{\lambda}, \bm{\mu}) \otimes L (\bm{\chi})$ is a simple quotient of $\widehat{R}(\psi)$, where $\psi$ is a Whittaker function determined by \eqref{uvjeti-2}.
 \end{Theorem}

\begin{proof} First, we note that theorem was already proved in \cite[Theorem 10.4]{ALZ-16} in the special case $\bm{\lambda} = (\lambda_0, 0, 0, \dots )$, $\bm{\mu} = (\mu_1, 0, 0, \dots )$, $\lambda_0 \ne 0$ and $\bm{\chi}$ arbitrary. Our situation is slightly more general, but the same proof applies here. Let us just present the basic idea for the reader’s convenience.

 Let $ h_1 =-2 {:} a a^*{:}$, ${\mathfrak b}_1 = {\C} e + {\C} h_1$ be a two dimensional Borel Lie algebra with bracket
$[h_1, e] = 2 e$ and
$\widehat {{\mathfrak b}} _1 = {\mathfrak b}_1 \otimes {\C}[t, t^{-1}] + {\C} K$ its affinization. Using the same arguments as in the proof of \cite[Proposition~10.2]{ALZ-16} we get
\begin{itemize}\itemsep=0pt
\item $M_1(\bm{\lambda}, \bm{\mu})$ is an irreducible $\widehat {{\mathfrak b}} _1$-module.
\end{itemize}

Next, we consider the Borel Lie algebra $\mathfrak b=\mathbb C e+\mathbb C h$, where
$h=h_1+b$,
with $b$ the Heisenberg field generator from Section~\ref{section-Wakimoto-1}, and its affinization $\widehat{ \mathfrak b} = {\mathfrak b} \otimes {\C} [t,t^{-1}] + {\C} K \subset \widehat {\mathfrak{sl}}_2$. Then we consider $M_1 (\bm{\lambda}, \bm{\mu}) \otimes L (\bm{\chi})$ as an \smash{$\widehat{ \mathfrak b}$}-module. The irreducibility of $M_1(\bm{\lambda}, \bm{\mu})$ as a~$\widehat {{\mathfrak b}} _1$\nobreakdash-mod\-ule is not affected if we twist the action of $h_1(z)$ by $\bm{\chi}(z)$.
Thus, we get
\begin{itemize}\itemsep=0pt
\item $ M_1(\bm{\lambda}, \bm{\mu}) \otimes L(\bm{\chi})$ is an irreducible $\widehat{ \mathfrak b}$-module.
\end{itemize}
 Since $\widehat{ \mathfrak b} $ is a Lie subalgebra of $ \widehat {\mathfrak{sl}}_2$ we conclude that $M_1(\bm{\lambda}, \bm{\mu}) \otimes L(\bm{\chi})$ is irreducible $ \widehat {\mathfrak{sl}}_2$-module and therefore irreducible $V^{\kappa} (\mathfrak{sl}_2)$-module. The proof follows.
 \end{proof}

\subsection*{Acknowledgements}

 We would like to thank the referees
for their valuable comments.
The authors are partially supported by the Croatian Science Foundation under the project IP-2022-10-9006 and by the project ``Implementation of
cutting-edge research and its application as part of the Scientific
 Center of Excellence for Quantum and Complex Systems, and
 Representations of Lie Algebras'', Grant No. PK.1.1.10.0004,
 co-financed by the European Union through the European Regional
 Development Fund -- Competitiveness and Cohesion Programme 2021--2027.
 This paper was also supported in part by the European Union -- NextGenerationEU through the National Recovery and Resilience Plan 2021--2026. Institutional grant of University of Zagreb Faculty of Science Strengthening scientific production, international presence and social impact of mathematical research (IK IA 1.1.3. Impact4Math).

\pdfbookmark[1]{References}{ref}
\LastPageEnding

\end{document}